\documentclass[11pt]{llncs}
\usepackage{graphicx}
\usepackage{cite}
\usepackage{url}
\usepackage{times}
\usepackage{mathptm}

\DeclareSymbolFont{AMSb}{U}{msb}{m}{n}
\DeclareSymbolFontAlphabet{\Bbb}{AMSb}
\def\R{\ensuremath{\Bbb R}}

\let\oldendproof\endproof
\def\endproof{\qed\oldendproof}

\begin{document}

\title{Cubic Partial Cubes from Simplicial Arrangements} 

\author{David Eppstein}

\institute{Computer Science Department\\
Donald Bren School of Information \& Computer Sciences\\
University of California, Irvine\\
\email{eppstein@uci.edu}}

\maketitle   

\begin{abstract}
We show how to construct a cubic partial cube from any simplicial arrangement of lines or pseudolines in the projective plane.  As a consequence, we find nine new infinite families of cubic partial cubes as well as many sporadic examples.
\end{abstract}

\section{Introduction}

\begin{figure}
\centering
\includegraphics[width=2.5in]{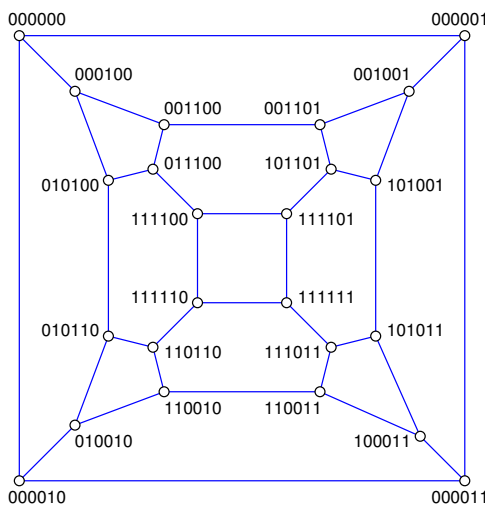}
\caption{A cubic partial cube, with vertex labels.}
\label{fig:p4}
\end{figure}

\begin{figure}[t]
\centering
\includegraphics[width=4in]{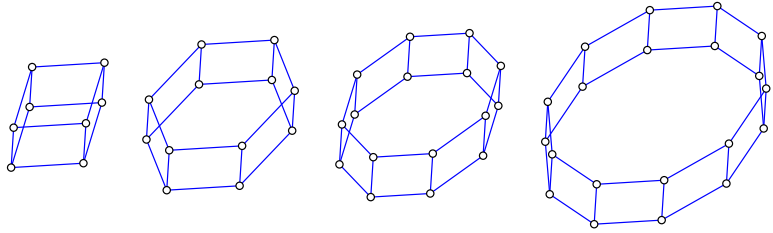}
\caption{The known infinite family of cubic partial cubes: prisms over even polygons.}
\label{fig:prisms}
\end{figure}

A {\em partial cube}~\cite{GraPol-BSTJ-71,Djo-JCTB-73} or {\em binary Hamming graph} is an undirected  graph, the vertices of which can be labeled by binary vectors in such a way that the distance between any two vertices in the graph is equal to the Hamming distance between the two vertices' labels.  In other words, the graph can be embedded isometrically onto a hypercube.  A partial cube, together with an isometric labeling of its vertices, is illustrated in Figure~\ref{fig:p4}. There has been much study of these graphs (e.g., \cite{AurHag-MST-95,DezLau-97,Epp-EJC-05,ImrKla-00,ImrKla-EJC-96,Win-DAM-84}) and their combinatorial enumeration~\cite{Wei-DM-92,BonKlaLip-AJC-03,BreKlaLip-EuJC-04,KlaLip-DM-03,KlaShp-ms-05,math.CO/0411359}; an interesting question in this area concerns classifying all {\em cubic} (that is, 3-regular) partial cubes~\cite{BonKlaLip-AJC-03,BreKlaLip-EuJC-04,KlaLip-DM-03,KlaShp-ms-05,math.CO/0411359}.  As Klav{\v z}ar and Shpectorov~\cite{KlaShp-ms-05} note, there are very large numbers of partial cubes (the subclass of median graphs corresponds roughly in numbers to the triangle-free graphs~\cite{ImrKlaMul-SJDM-99}) so the difficulty of finding cubic partial cubes comes as somewhat of a surprise.  The known cubic partial cubes fall into one infinite family (the prisms over even polygons, shown in Figure~\ref{fig:prisms}), together with 36 sporadic examples not known to belong to any infinite family~\cite{KlaShp-ms-05}.

\begin{figure}[t]
\centering
\includegraphics[width=2.5in]{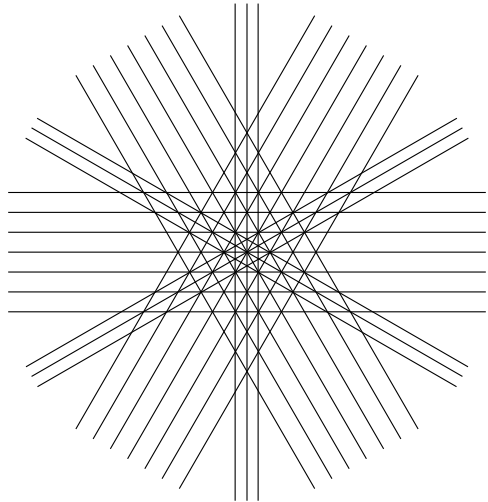}
\caption{A simplicial line arrangement, including the line at infinity.}
\label{fig:737373}
\end{figure}

There are many known similar situations in combinatorics, of objects that can be enumerated as one or a few infinite families together with finitely many sporadic examples.  Another such is that of {\em simplicial arrangements} of lines in the real projective plane; that is, finite sets of lines such that each of the cells of the arrangement (including the cells meeting or passing through the line at infinity) is a triangle. One such arrangement is shown in Figure~\ref{fig:737373}; note that, along with the many finite triangles in the figure, each half-strip cell bounded by a pair of parallel lines is projectively a triangle with a vertex at infinity, while the twelve wedge-shaped cells form triangles having the line at infinity along one of their sides.  It is necessary to include the line at infinity in this figure, for otherwise opposite pairs of wedge-shaped cells would form quadrilaterals in the projective plane.  The enumeration of simplicial arrangements has been well studied~\cite{Gru-72}, and three infinite families of such arrangements are known, together with 91 sporadic examples.

There is a standard construction that forms a dual partial cube from any line or hyperplane arrangement~\cite{Ovc-DAM-?}, by forming a vertex for each cell and connecting two vertices by an edge whenever the corresponding two cells border each other along a hyperplane.  However, this contruction applies to arrangements in affine spaces, not projective ones.  When applied to a simplicial line arrangement, it produces a partial cube with degree three vertices for each finite triangle, however some cells (e.g. the wedge to the left of the leftmost vertex of the arrangement) will lead to degree two vertices.

In this paper we overcome this difficulty of relating simplicial arrangements to cubic partial cubes.  We show how to connect two copies of the partial cube corresponding to the affine part of a simplicial arrangement to each other, to provide a construction of a cubic partial cube from any simplicial arrangement.  Equivalently, our construction can be seen as forming the dual to a simplicial central arrangement of planes in $\R^3$, or the dual to a simplicial line arrangement in the {\em oriented projective plane}~\cite{Sto-91}.  Using this construction, we form many new cubic partial cubes.  One of the three known infinite families of simplicial arrangements (the near-pencils) leads to a known set of cubic partial cubes (the prisms) but the other two lead to two new infinite families of cubic partial cubes.  In addition, the sporadic simplicial arrangements lead to many new sporadic cubic partial cubes.

We then generalize our construction to simplicial {\em pseudoline} arrangements.
Seven additional infinite families of such arrangements are known, together with many more sporadic examples; again, each leads to a new cubic partial cube or infinite family of cubic partial cubes.

Finally, we discuss cubic partial cubes that do not come from arrangements in this way.
We describe a construction for gluing together the duals to two different affine arrangements, or of two rotated copies of the same arrangement; our construction leads to several new cubic partial cubes that are not dual to arrangements.

\section{Preliminaries}

The {\em Djokovi{\'c} relation} is a binary relation on the edges of a graph, under which two edges $(v,w)$ and $(x,y)$ are related if and only if $d(v,x)=d(w,y)=d(v,y)-1=d(w,x)-1$ for some ordering of the endpoints of the two edges.
We recall~\cite{Djo-JCTB-73} that a graph is a partial cube if and only if its Djokovi{\'c} relation is an equivalence relation; in that case we call its equivalence classes {\em Djokovi{\'c} classes}.
Each Djokovi{\'c} class forms a cut separating the graph into two connected components.
If a partial cube is labeled by binary vectors in such a way that graph distance equals Hamming distance, then the endpoints of each edge have labels that differ in a single coordinate, and we can group the edges into classes according to which coordinate their endpoint labels differ in; these classes are exactly the Djokovi{\'c} classes.  Conversely, from the Djokovi{\'c} classes we can form a binary labeling by assigning one coordinate per class, with a coordinate value of zero on one side of the cut formed by the class and a coordinate value of one on the other side of the cut.

\section{Partial Cubes from Line Arrangements}

We begin by recalling~\cite{Ovc-DAM-?} the construction of a partial cube from a hyperplane arrangement in $\R^d$.  Our construction of cubic partial cubes from simplicial line arrangements will use this construction in $\R^3$.

\begin{lemma}
\label{lem:affine-pcube}
Let $A$ be a finite set of hyperplanes in $\R^d$, and form a graph $G_A$ that has one vertex for each $d$-dimensional cell of the arrangement, and that connects two vertices by an edge whenever the corresponding two cells meet along a $(d-1)$-dimensional face of the arrangement.  Then $G_A$ is a partial cube.
\end{lemma}

\begin{proof}
For each halfplane $H_i$ of $A$, choose $P_i$ arbitrarily to be one of the two halfspaces bounded by $H_i$.  For each cell $c$ of the arrangement, label the corresponding vertex by the binary vector $b_0b_1b_2\ldots$, where $b_i=1$ if $c$ is contained in $P_i$ and $b_i=0$ otherwise.

If two vertices of $G_A$ are connected by an edge, the corresponding cells meet along a face belonging to a hyperplane $H_i$, and exactly one of the two vertices has $b_i=1$; for all other $H_j$, $j\ne i$, the two cells are both on the same side of $H_j$.  Therefore, the Hamming distance between the labels of two adjacent vertices is exactly one, and the Hamming distance between any two vertices is at most equal to their graph distance.

Conversely, suppose that $v$ and $v'$ are two vertices of $G_A$ corresponding to cells $c$ and $c'$ of $A$.  Choose $p$ and $p'$ interior to these two cells such that all points of the line connecting them belong to at most one hyperplane of $A$; this can be done as, for any $p$, the set of $p'$ not satisfying this condition forms a subset of measure zero of $c'$.  Form a path from $v$ to $v'$ in $G_A$ corresponding to the sequence of cells crossed by line segment $pp'$.  This line segment crosses only hyperplanes $H_i$ of $A$ for which the labels $b_i$ at $v$ and $v'$ differ, so this path has length exactly equal to the Hamming distance of the two labels.

We have shown that, for every two vertices in $G_A$, the Hamming distance of their labels equals the length of the shortest path connecting them, so $G_A$ is a partial cube.
\end{proof}

\begin{theorem}
\label{thm:line-pcube}
Let $A$ be a simplicial line arrangement in the projective plane.
Then there corresponds to $A$ a cubic partial cube $C_A$,
with twice as many vertices as $A$ has triangles.
\end{theorem}

\begin{proof}
Embed the projective plane as the plane $z=1$ in $R^d$, and form an arrangement
$\hat A$ that has, for each line $\ell$ of the arrangement $A$, a plane through $\ell$ and the origin.  If $A$ contains the line at infinity, add correspondingly to $\hat A$ the plane $z=0$.  Then $\hat A$ is an arrangement of planes, in which each cell is an infinite triangular cone.
The graph $C_A=G_{\hat A}$ is, by Lemma~\ref{lem:affine-pcube}, a cubic partial cube.
\end{proof}

\begin{theorem}
If $A$ is a simplicial line arrangement, then the graph $C_A$ constructed in Theorem~\ref{thm:line-pcube} is planar.
\end{theorem}

\begin{proof}
If we intersect the arrangement $\hat A\subset\R^3$ with the unit sphere, we obtain an arrangement of great circles on the sphere with spherical-triangle faces.  The vertices and arcs of the arrangement on the sphere form a planar graph, and our construction $C_A$ is just the planar dual of this graph.  Thus, all graphs $C_A$ constructed via Theorem~\ref{thm:line-pcube} are planar.
\end{proof}

The geometry of points and great circles on the sphere can also be interpreted as a model of {\em oriented projective geometry}~\cite{Sto-91}, which consists of two signed points for each point in the more standard unoriented projective plane.  We will use this interpretation later when we generalize from lines to pseudolines.

\begin{lemma}
\label{lem:unique-embed}
If $A$ is an arrangement of lines that is not a pencil, in the oriented projective plane,
then the planar graph of vertices and edges of $A$ has a unique planar embedding.
\end{lemma}

\begin{proof}
Consider any face $f$ of $A$.  Then $f$ is {\em nonseparating}: that is, every two edges of
$A\setminus f$ can be connected by a path that does have as its interior vertices any vertices of $f$.  For, if one has a path connecting two edges, passing through $f$, the portion of the path passing through edges of $f$ can be replaced by paths through the remaining portions of the great circles containing those edges.

In any embedding of any planar graph, any nonseparating cycle must be a face.  Thus, all faces of the arrangement must be faces in any planar embedding of the same graph.  But, by Euler's formula, the number of faces in all planar embeddings must be equal.  Therefore, any embedding of the planar graph of $A$ has the same set of cycles as its set of faces, and is therefore combinatorially equivalent to the arrangement's embedding.
\end{proof}

\begin{theorem}
If $A$ and $A'$ are simplicial line arrangements, then the
cubic partial cubes $C_A$ and $C_{A'}$ are isomorphic graphs if and only if $A$ and $A'$ are combinatorially equivalent as arrangements.
\end{theorem}

\begin{proof}
Our construction depends only on the combinatorial type of the arrangement, not on the geometric positioning of its lines, so if $A$ and $A'$ are combinatorially equivalent then $C_A$ and $C_{A'}$ are isomorphic.

In the other direction, suppose that $C_A$ and $C_{A'}$ are isomorphic as graphs.
By Lemma~\ref{lem:unique-embed}, $A$ and $A'$, and therefore also their duals $C_A$ and $C_{A'}$, are uniquely embeddable, and their embeddings must be combinatorially equivalent.
Therefore the planar duals of $C_A$ and $C_{A'}$ are also isomorphic to each other,
with combinatorially equivalent embeddings.
But these duals are just the graphs of the two arrangements $A$ and $A'$, and the intersection patterns of the lines of the arrangements can be recovered from curves in these graphs that, at each vertex, pass through oppositely situated edges in the embedding.  Therefore also $A$ and $A'$ are combinatorially equivalent.
\end{proof}

Thus, we can construct a distinct cubic partial cube for each possible simplicial line arrangement.

\section{The Three Infinite Families}

\begin{figure}[t]
\centering
\includegraphics[width=1.85in]{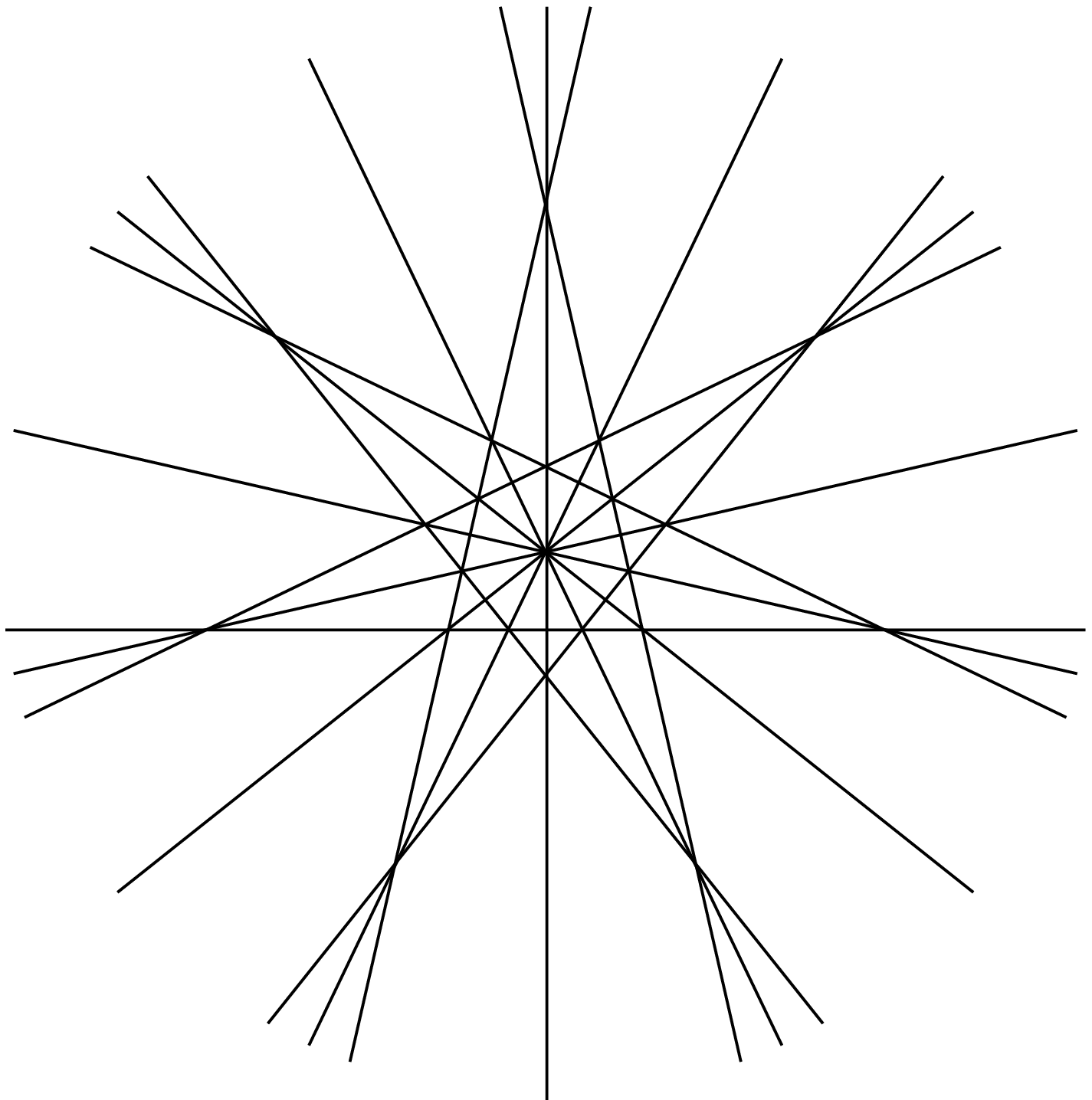}\hspace{0.25in}
\includegraphics[width=1.85in]{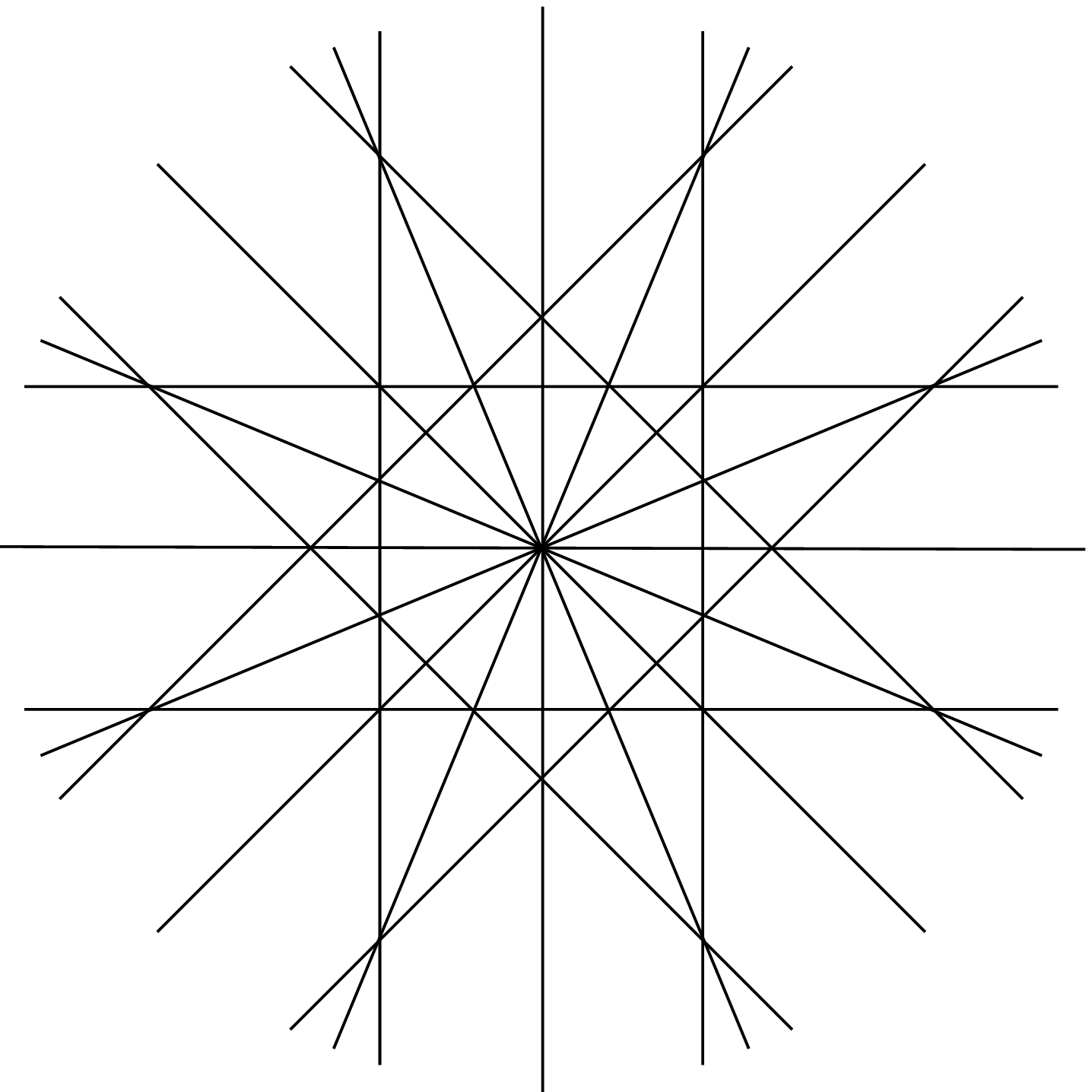}\hspace{0.25in}
\includegraphics[width=1.85in]{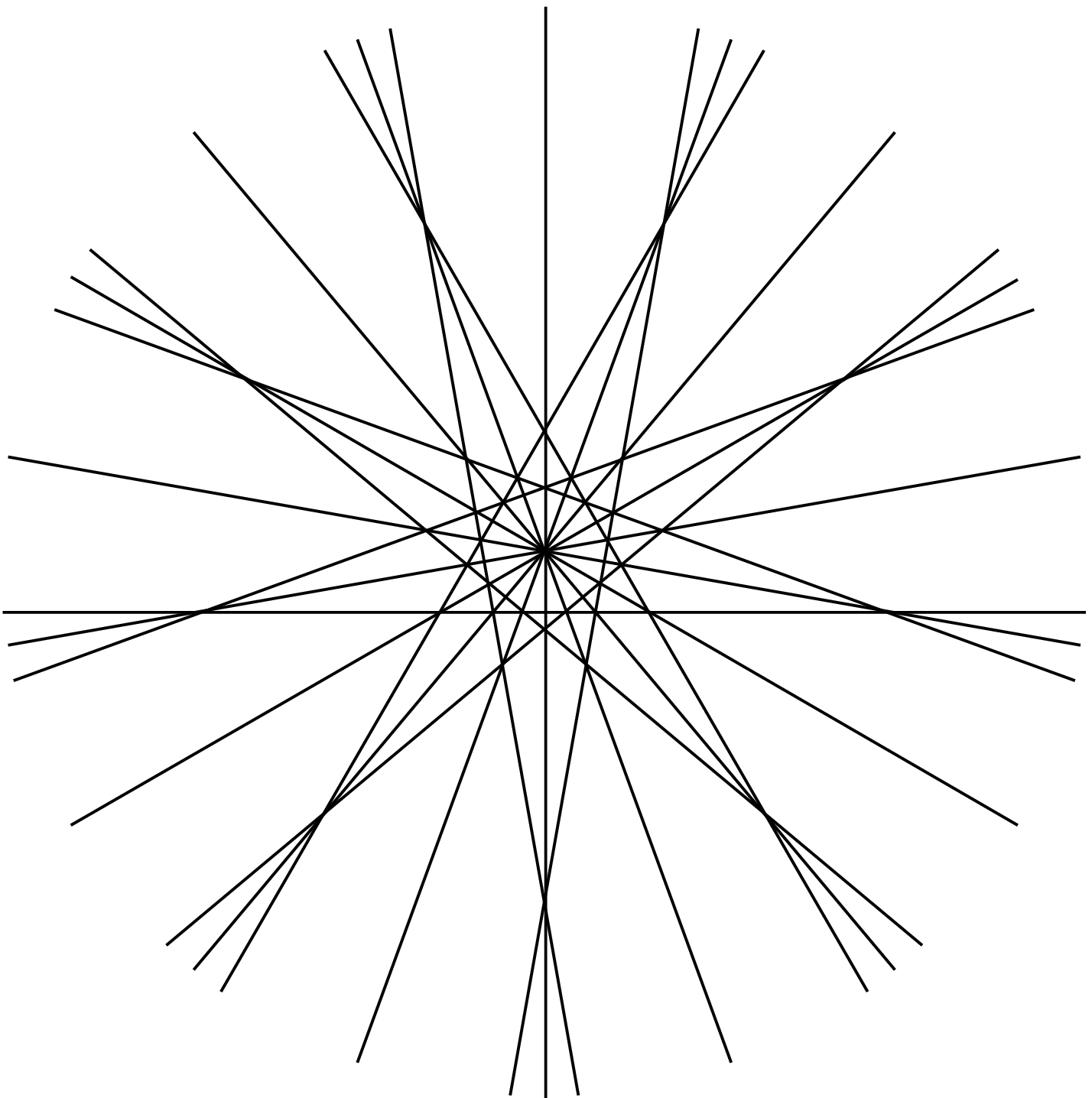}

\vspace{0.25in}

\includegraphics[width=1.85in]{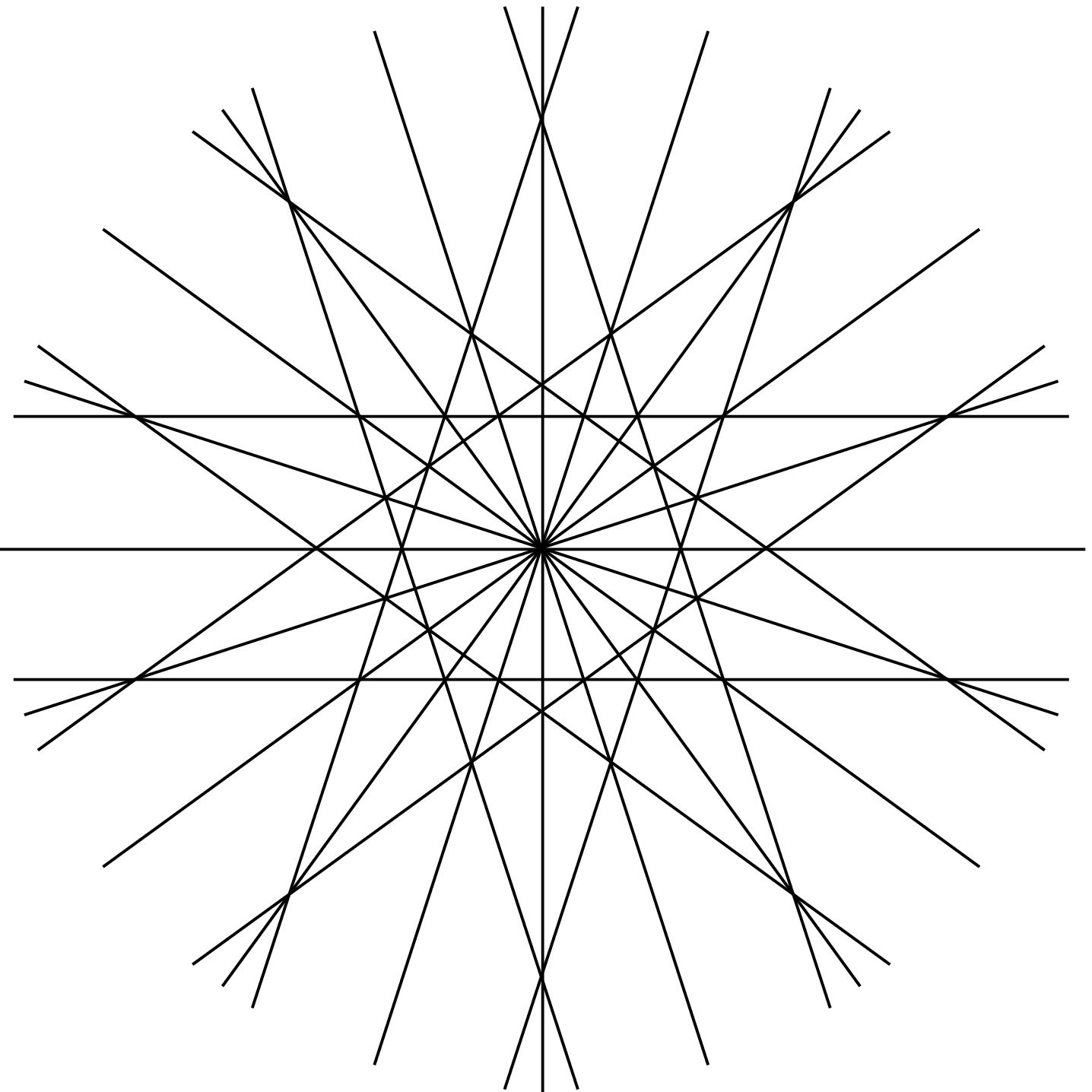}\hspace{0.25in}
\includegraphics[width=1.85in]{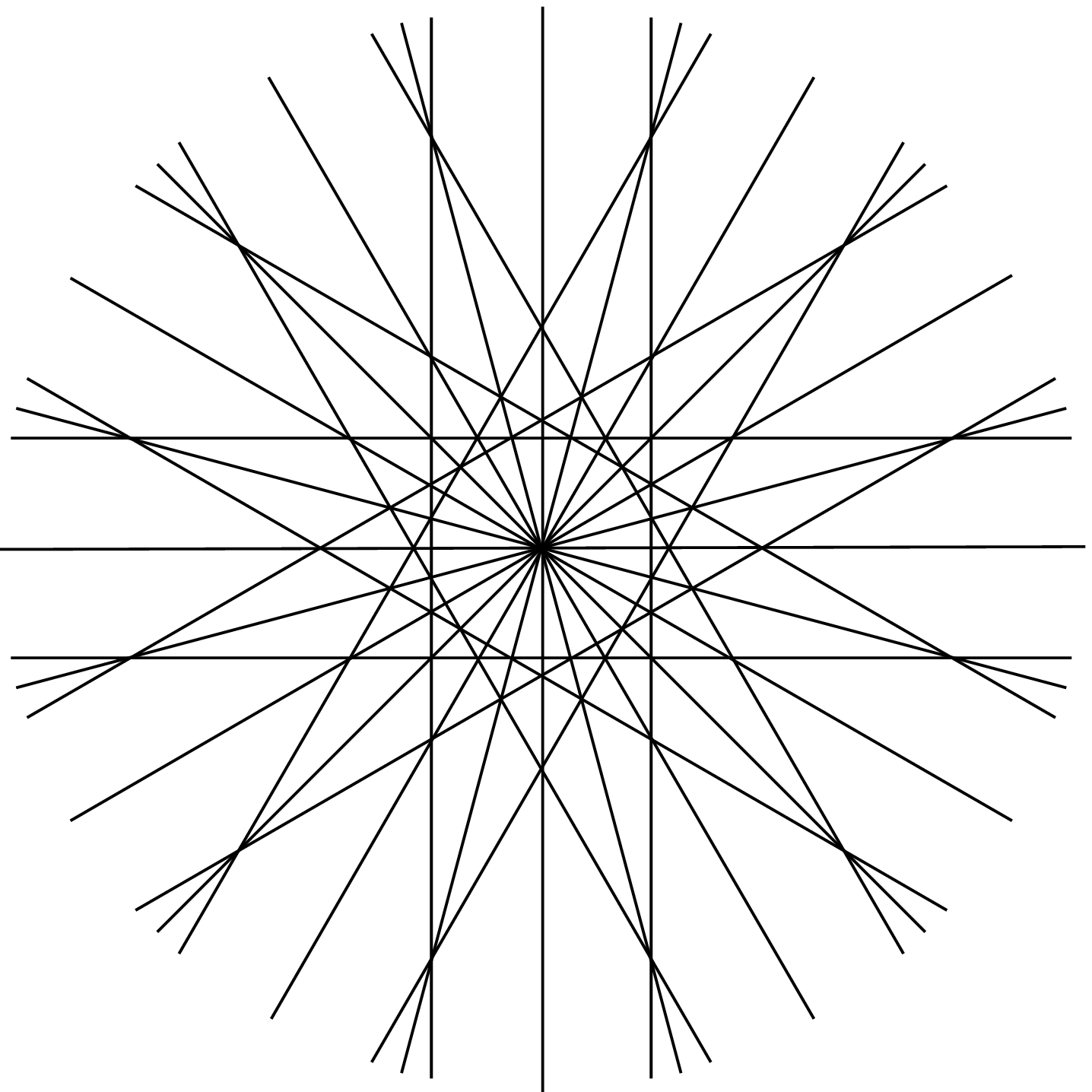}\hspace{0.25in}
\includegraphics[width=1.85in]{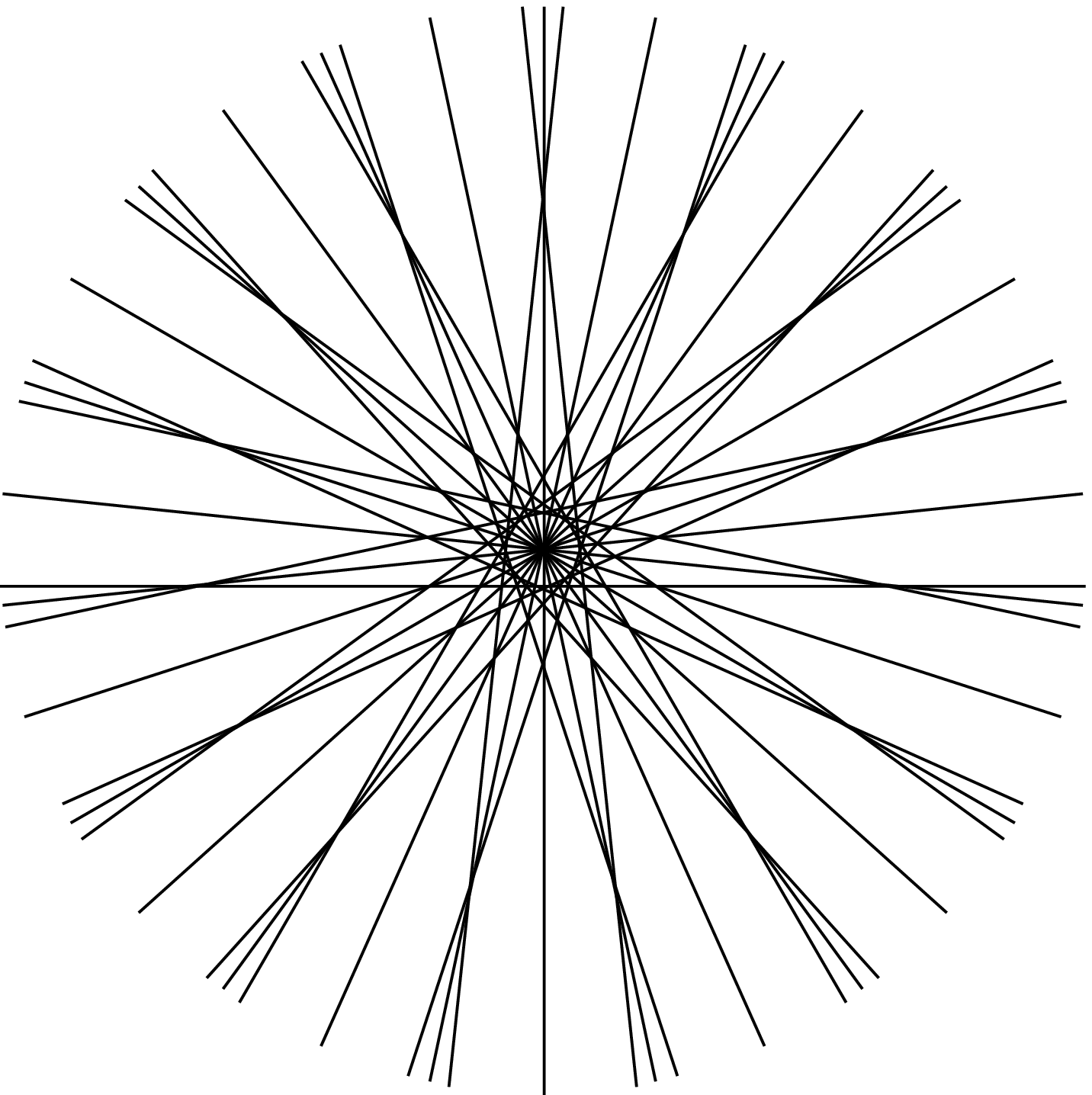}
\caption{The simplicial arrangements $R(14)$, $R(16)$, $R(18)$,
$R(20)$, $R(24)$, and $R(30)$.}
\label{fig:R}
\end{figure}

There are three known infinite families of simplicial line arrangements~\cite{Gru-72}. The first are the {\em near-pencils}; the near-pencil of $n$ lines consists of $n-1$ lines through a single point (a {\em pencil}) together with an additional line not through that point.  We may form this arrangement most symmetrically by spreading the lines of the pencil at equal angles, and placing the additional line at infinity.  If $A$ is a near-pencil, the dual cubic partial cube $C_A$ constructed via Theorem~\ref{thm:line-pcube} is a known cubic partial cube, a prism (Figure~\ref{fig:prisms}).

The second infinite family is denoted $R(2k)$, and consists of $k$ lines formed by extending the sides of a regular $k$-gon, together with an additional $k$ lines formed by the $k$ axes of symmetry of the $k$-gon. When $k$ is odd, each of these axes passes through one vertex and one edge midpoint of the $k$-gon; when $k$ is even, they pass through two vertices or two edge midpoints.
Several members of this family are depicted in Figure~\ref{fig:R}.
$R(4k)$ has $2k(2k+1)$ triangles in the projective plane, and therefore leads to the construction of cubic partial cubes with $4k(2k+1)$ vertices and $4k$ Djokovi{\'c} classes of edges.
$R(4k+2)$ has  $2(k+1)(2k+1)$ triangles in the projective plane, and therefore leads to the construction of cubic partial cubes with $4(k+1)(2k+1)$ vertices and $4k+2$ Djokovi{\'c} classes of edges.

The third infinite family is denoted $R(4k+1)$, and consists of the $4k$ lines of $R(4k)$, together with a single additional line at infinity; thus, the arrangements $R(17)$, $R(21)$, and $R(25)$ can be formed by adding a line at infinity to the figures depicting $R(16)$, $R(20)$, and $R(24)$.
Adding a line at infinity to $R(4k+2)$ produces an arrangement that is not simplicial.
The simplicial arrangements $R(4k+1)$ have $4k(k+1)$ triangles in the projective plane,
and therefore lead to the construction of cubic partial cubes with $8k(k+1)$ vertices and $4k+1$ Djokovi{\'c} classes of edges.

\section{Pseudoline Arrangements}

\begin{figure}[t]
\centering
\includegraphics[width=2.5in]{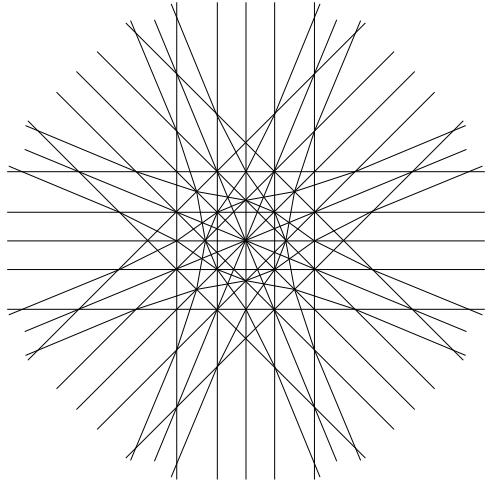}
\caption{A simplicial pseudoline arrangement (with the line at infinity).}
\label{fig:pseudo}
\end{figure}

Although there are many simplicial line arrangements known, there are even more simplicial {\em pseudoline arrangements}.  A pseudoline arrangement is a collection of curves in the projective plane, each topologically equivalent to a line (that is, closed non-self-crossing non-contractible curves), such that any two curves in the collection have a single crossing point.  Using pseudolines, we can form simplicial arrangements (that is, arrangements in which each cell is bounded by sides belonging to three curves) that may not be realizable as line arrangements, such as the one shown in Figure~\ref{fig:pseudo}.

The interpretation of Theorem~\ref{thm:line-pcube} as planar duality in the oriented projective plane generalizes to simplicial pseudoline arrangements, and produces a cubic graph $C_A$ from any such arrangement.
That is, we lift the arrangement from the projective plane to the oriented projective plane, and form the planar dual graph of the lifted arrangement.
When $A$ is a line arrangement, this coincides with the previous construction.

\begin{theorem}
\label{thm:pseudo-pcube}
The graph $C_A$ constructed as above from a pseudoline arrangement is a partial cube.
If the arrangement is simplicial, the partial cube is cubic.
\end{theorem}

\begin{proof}
Each pseudoline divides the oriented projective plane into two halfspaces. We may assign a zero to one of these halfspaces and a one to the other arbitrarily, and label each cell by the sequence of bits from the halfspaces containing it.  These labels differ by one across each edge of the arrangement, so the Hamming distance in $C_A$ is at least the graph distance.
It remains to show that any two vertices $v$ and $w$ in the graph have a path as short as the Hamming distance between their labels.

We extend the arrangement $A$ in the (unoriented) projective plane to a {\em spread}; that is, an infinite collection of pseudolines, each pair having a single crossing point as before, that covers the plane in the sense that every two points have a unique pseudoline connecting them.  It is known~\cite{GooPolWen-Comb-94} that every arrangement can be extended to a spread in this way.  We lift this spread to the oriented projective plane, and from the cells in $A$ in the oriented projective plane corresponding to $v$ and $w$ in the graph, we choose representative points $p_v$ and $p_w$ in general position, meaning that $p_v$ and $p_w$ are not both the image of a single point in the projective plane and that the pseudoline connecting these two points in the spread does not pass through any vertices of the arrangement.  We then find a path from $v$ to $w$ by following the sequence of cells crossed by the pseudoline segment from $p_v$ to $p_w$.  This pseudoline segment crosses each pseudoline of $A$ at most once, and if it does cross a pseudoline then that pseudoline's coordinate in the vertex labels will differ from before the crossing to after it. Therefore, the total number of steps can at most equal the number of coordinates on which the labels of $v$ and $w$ differ.
\end{proof}

\begin{figure}[t]
\centering
\includegraphics[width=3.5in]{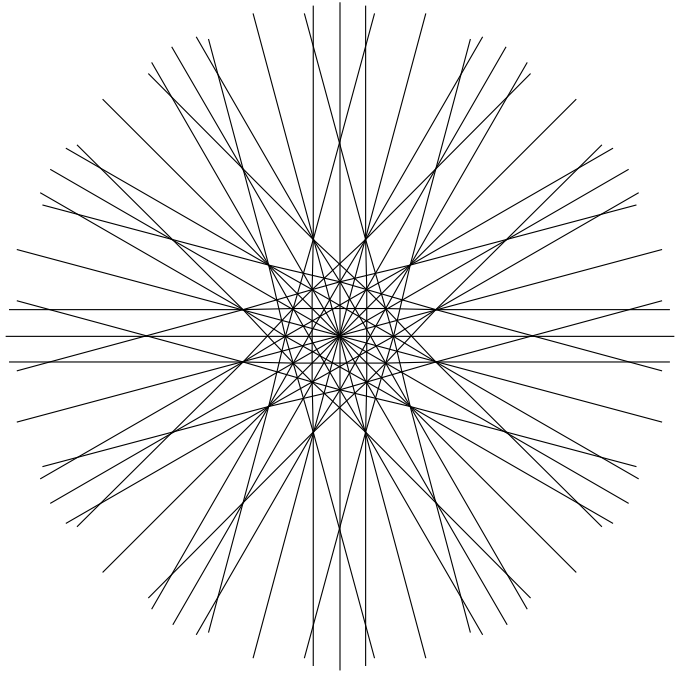}
\caption{An arrangement (including the line at infinity) belonging to an infinite family of pseudoline arrangements with $3k+1$ lines, for $k=4$ or $6$ mod $10$. This arrangement has 37 lines, and can be formed by interleaving two differently-scaled copies of $R(25)$. It can be realized with straight lines, but other members of the family cannot. All finite pseudolines in all members of this family occur in groups of three parallel pseudolines.}
\label{fig:P37}
\end{figure}

Gr\"unbaum~\cite{Gru-72} notes without reference or detail the existence of seven infinite families of simplicial pseudoline arrangements.  For instance, several such families may be formed by interleaving two differently-scaled copies of the line arrangements $R(n)$, as in Figure~\ref{fig:P37} and
Gr\"unbaum's figures 3.15-3.17.
Thus, together with the cubic partial cubes coming from line arrangements, we have constructed a total of ten infinite families of cubic partial cubes.  In addition, many additional sporadic examples of simplicial pseudoline arrangements are known.

\section{Geometric representation and non-arrangement-based partial cubes}

\begin{figure}[t]
\centering
\includegraphics[width=1.25in]{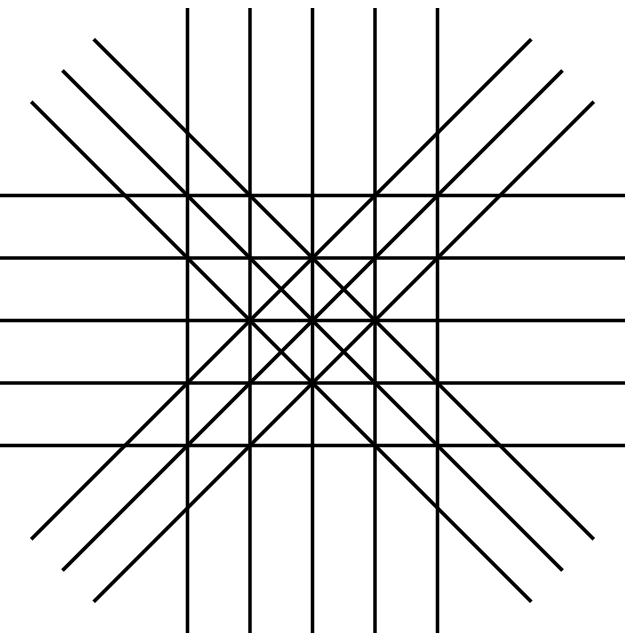}\hspace{0.25in}
\includegraphics[width=1.5in]{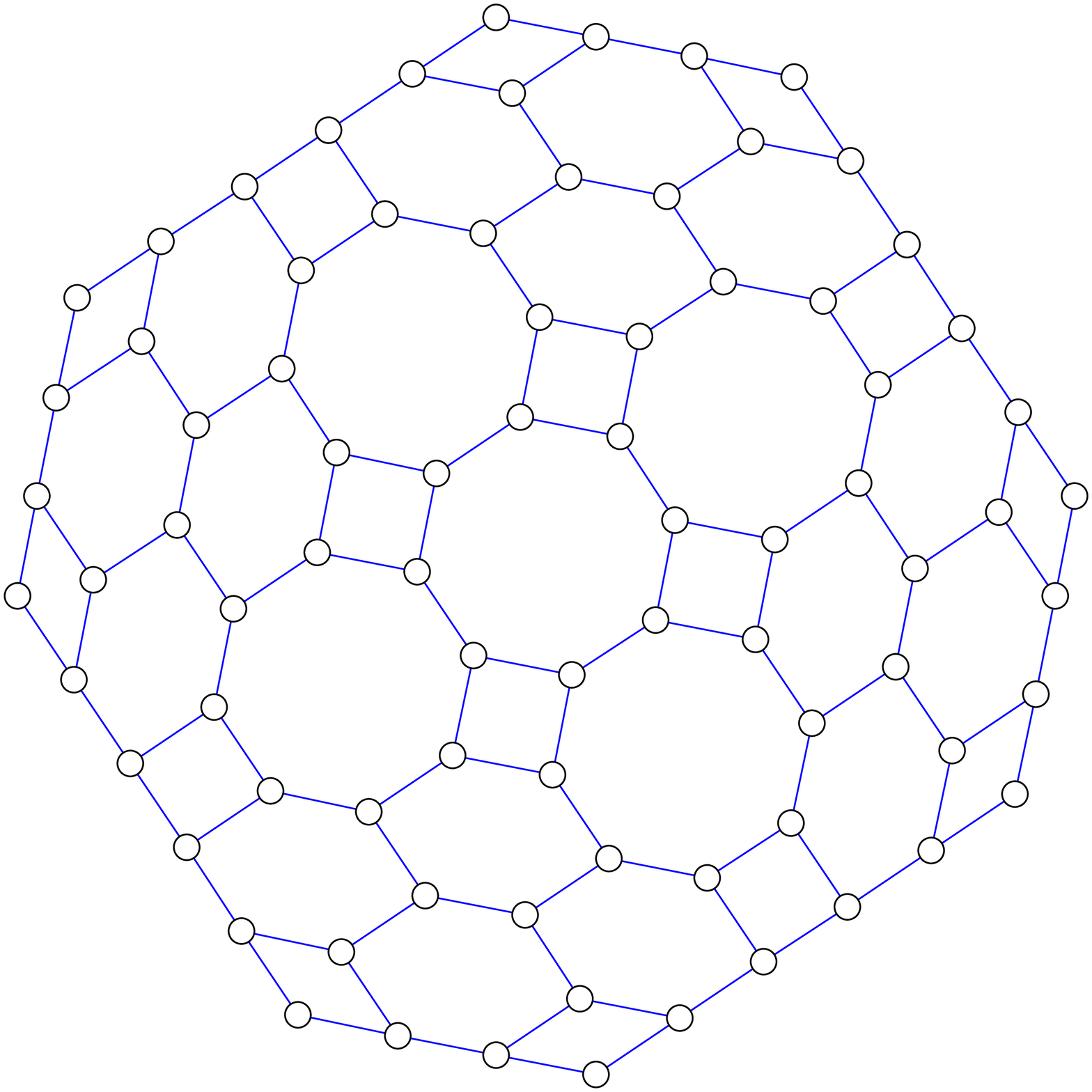}\hspace{0.25in}
\includegraphics[width=1.5in]{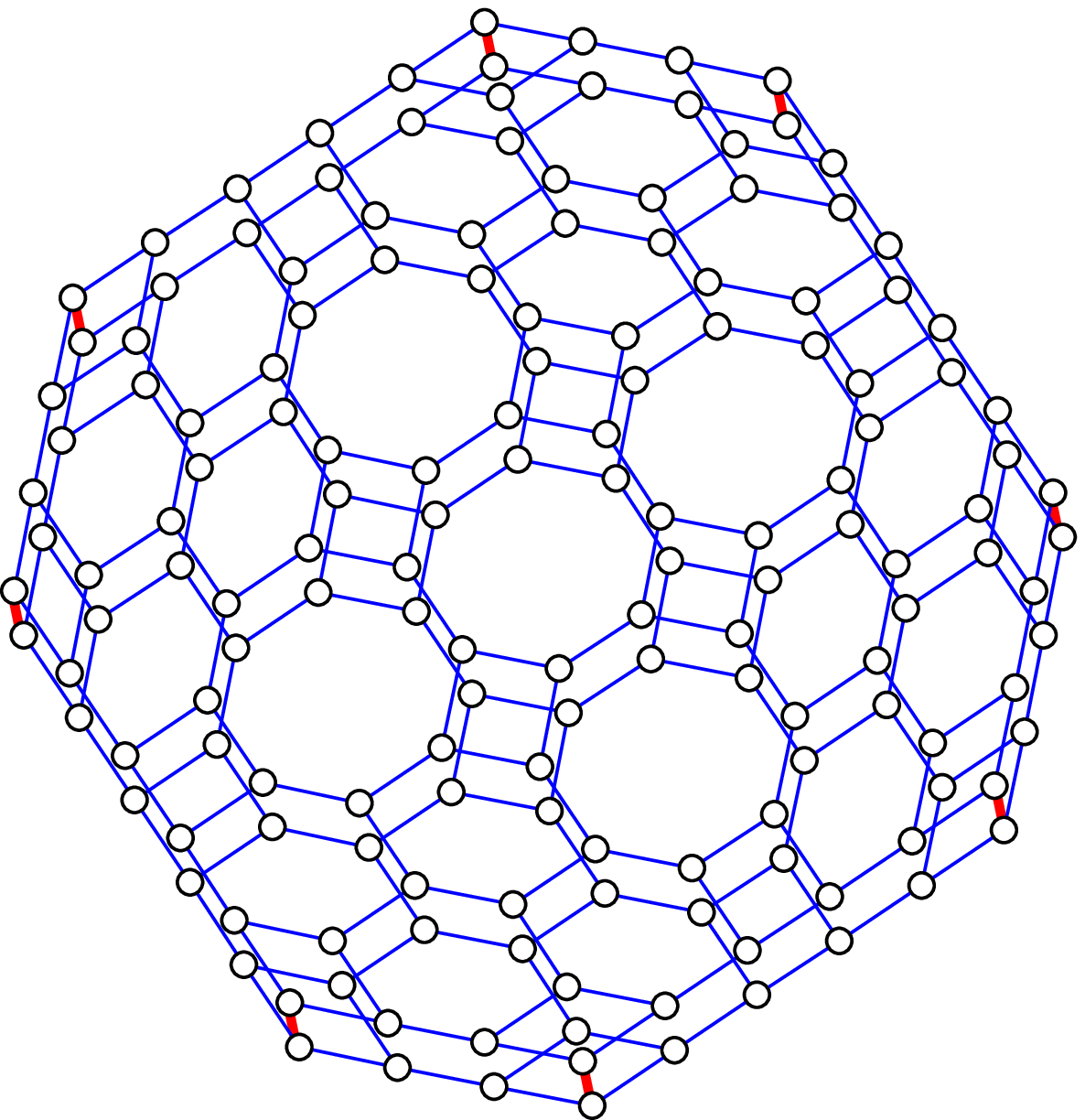}

\vspace{0.25in}

\includegraphics[width=1.25in]{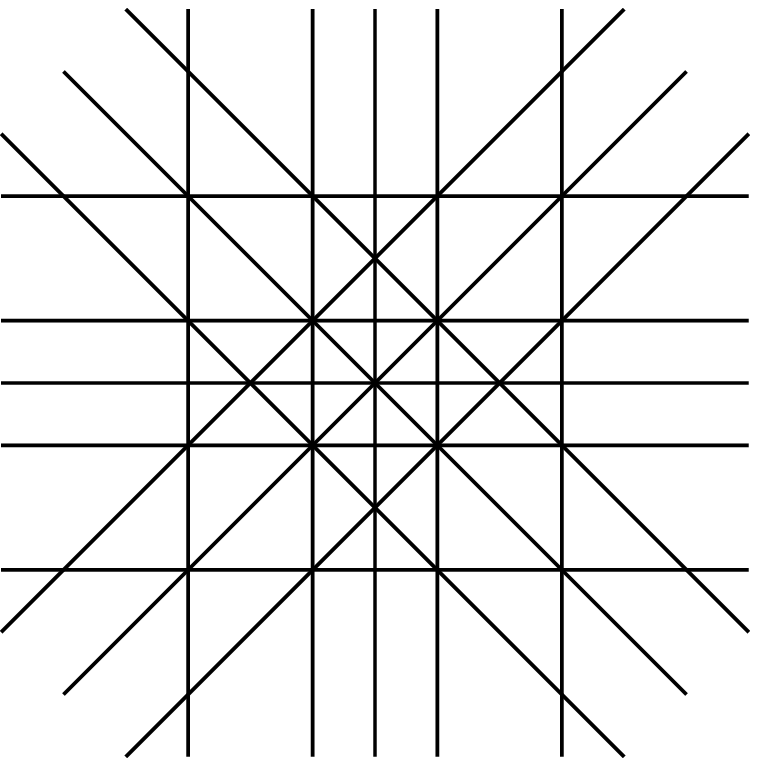}\hspace{0.25in}
\includegraphics[width=1.5in]{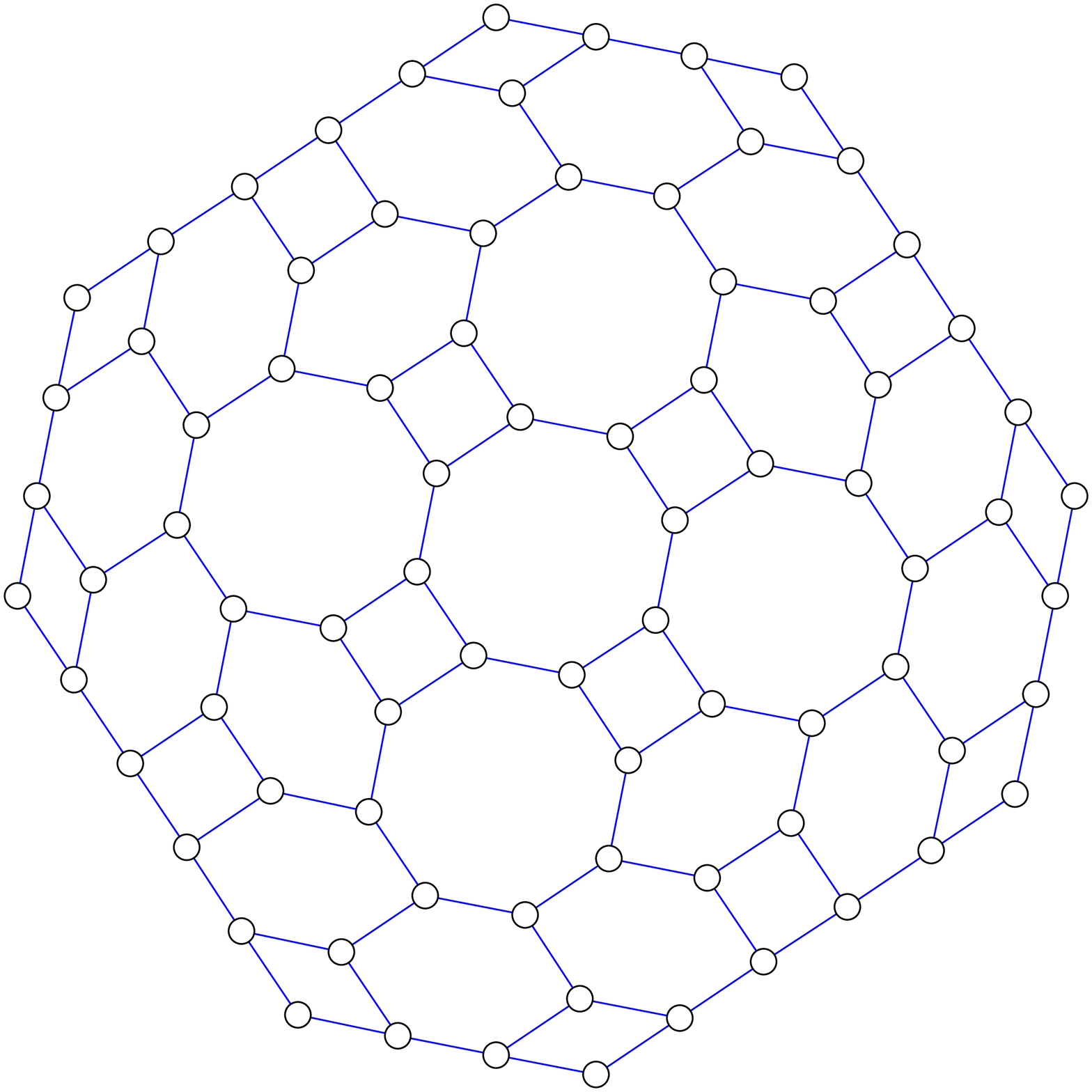}\hspace{0.25in}
\includegraphics[width=1.5in]{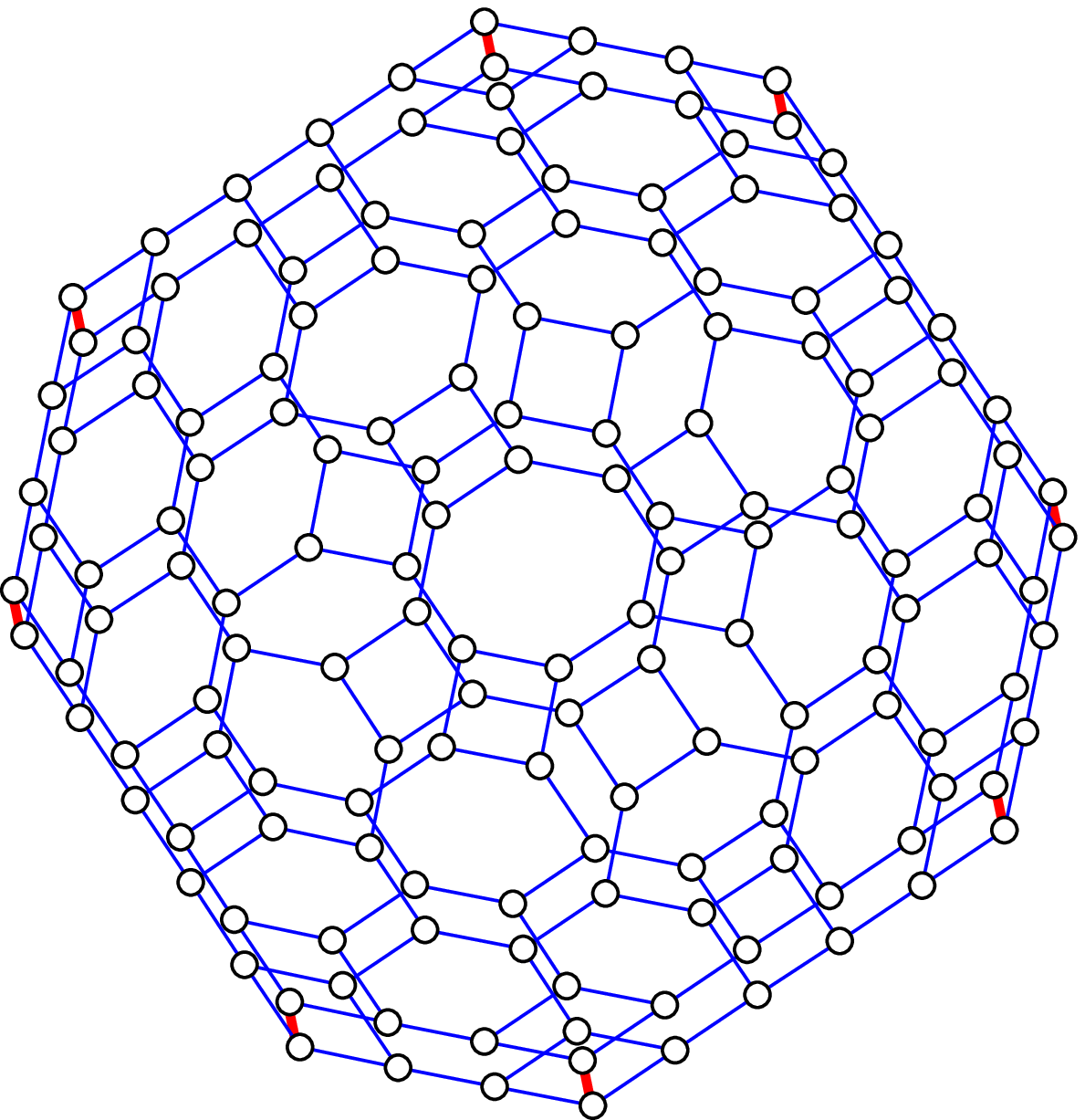}
\caption{Top: a simplicial arrangement (left), dual symmetric-faced planar drawing (center), and partial cube formed by connecting two copies of the drawing (right).
Bottom: a second simplicial arrangement with the same numbers of parallel lines of each slope (left), dual drawing (center), and non-arrangement-based partial cube formed by connecting the drawings of the first and second arrangements.}
\label{fig:AB}
\end{figure}

In an earlier paper~\cite{Epp-GD-04}, we described a method for finding planar drawings of partial cubes (when such drawings exist) in which all internal faces are centrally symmetric strictly convex polygons, based on a duality between these drawings and a generalization of pseudoline arrangements.  In this drawing method, as applied to the partial cubes dual to Euclidean line arrangements, we draw a vertex for each cell of the arrangement, and connect any two adjacent vertices by a unit length edge oriented perpendicularly to the line separating the two cells corresponding to the vertices.  The properties of partial cubes can be used to show that these rules define consistent placements for all vertices and edges of the graph; examples of drawings constructed in this way can be seen in the middle two parts of Figure~\ref{fig:AB} and in the right side of Figure~\ref{fig:Q}.
The resulting drawing is a so-called {\em zonotopal tiling} in which an outer centrally symmetric convex polygon is partitioned into smaller centrally symmetric convex polygons~\cite{FelWei-DAM-01}.

Unfortunately, this method is specific to Euclidean arrangements, although it works equally well for lines or pseudolines. It does not work for the oriented projective arrangements we use to construct cubic partial cubes.  No cubic partial cube can have a planar drawing in which all faces are centrally symmetric strictly convex polygons, because such a drawing can be shown to be dual to a certain kind of arrangement~\cite{Epp-GD-04}, and an arrangement cell to the left of the leftmost arrangement vertex would have to correspond to a vertex of degree two in the graph.

\begin{figure}[t]
\centering
\includegraphics[width=1.75in]{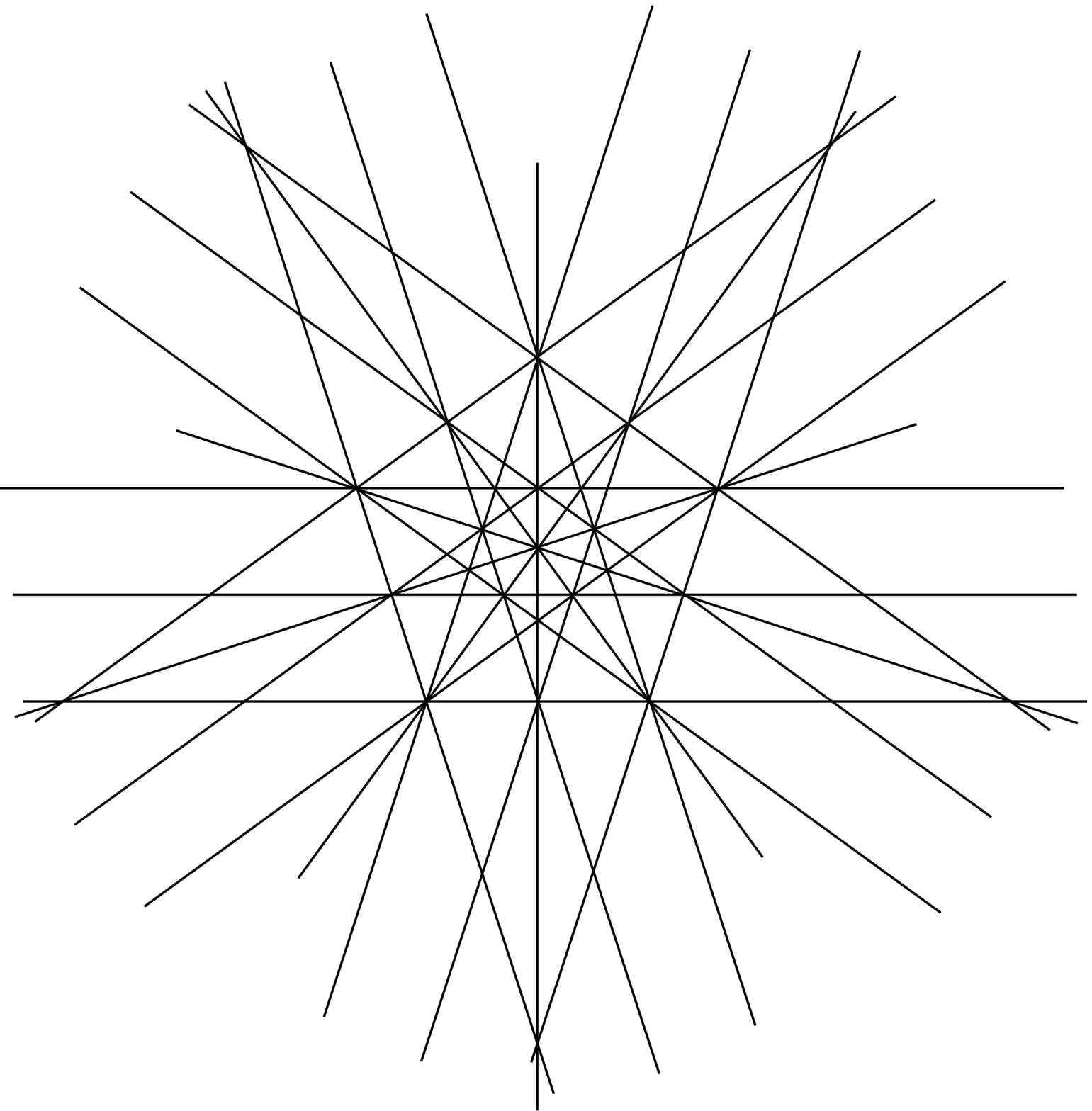}\hspace{0.25in}
\includegraphics[width=1.65in]{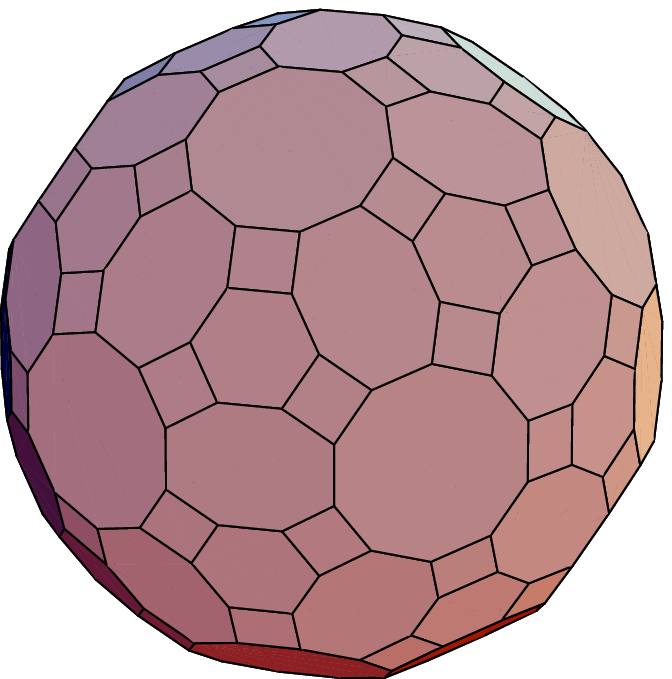}
\caption{A simplicial line arrangement (left) and the corresponding zonohedron (right, from~\cite{Epp-MER-96}).}
\label{fig:zono}
\end{figure}

On the other hand, when $A$ is a line arrangement, it is possible to construct a three-dimensional representation of $C_A$ with symmetric strictly convex faces: a {\em zonohedron}
(Minkowski sum of line segments)~\cite{Epp-MER-96}.
To do so, lift each line in the plane to a plane through the origin in $\R^3$, as in the proof of Theorem~\ref{thm:line-pcube}, and form the Minkowski sum of unit vectors perpendicular to each plane.
The resulting shape is a convex polyhedron, in which each face is a strictly convex and centrally symmetric polygon.  If one partitions the points of the unit sphere, according to which face of the polyhedron has the largest dot product with each point, the resulting partition forms an arrangement of great circles, equal to the great circle arrangement formed by intersecting the planes through the origin with the unit sphere.  Therefore, the zonohedron is the planar dual to the great circle arrangement, and its vertices and edges for a rpresentation of the cubic partial cube $C_A$.
However, this representation is only possible for $C_A$ when $A$ is a line arrangement, and not when $A$ is a pseudoline arrangement. Also, it requires prior knowledge of the coordinates of the lines in the arrangement, while our planar symmetric-faced drawing method needs only the graph structure of the underlying partial cube.

\begin{figure}[t]
\centering
\includegraphics[width=4.5in]{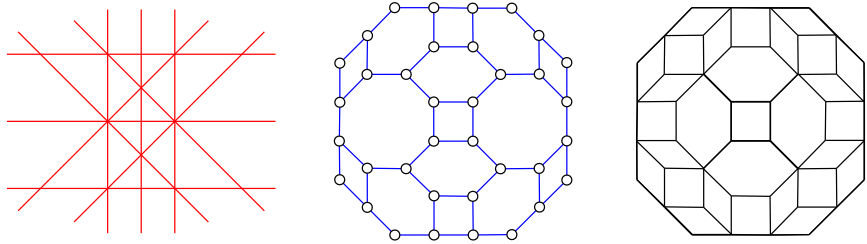}
\caption{A simplicial arrangement (left) and its dual drawing (center).  We can form a non-arrangement-based partial cube by connecting the drawing to a $90^\circ$ rotation of itself; at right is the zonotopal tiling formed by overlaying the drawing and its rotation..}
\label{fig:Q}
\end{figure}

We may combine these two methods to achieve three-dimensional representations with symmetric faces, most but not all of which are strictly convex, even for the graphs $C_A$ coming from pseudoline arrangements~$A$.  To do so, from a simplicial arrangement that includes the line at infinity, begin by forming the zonotopal tiling dual to the affine part of the arrangement, as in the center parts of Figure~\ref{fig:AB}.  This tiling covers a convex polygon in the plane, with degree three vertices except at the polygon's corners, which have degree two.  Place one copy of this tiling in $\R^3$, parallel to a centrally reflected copy of the tiling, and connect corresponding corners of the two copies (Figure~\ref{fig:AB}, right top).  The result is a three-dimensional representation of the graph, in which the planar faces corresponding to finite vertices of the arrangement are drawn as strictly convex polygons while the faces corresponding to vertices at infinity are drawn as subdivided rectangles.

In certain cases, we may apply a similar gluing approach to form cubic partial cubes that are not of the form $C_A$ for a single arrangement~$A$.  For instance, the top and bottom parts of Figure~\ref{fig:AB} depict two different line arrangements that have the same numbers of parallel lines of each slope (the left parts of the figure), so that their corresponding zonotopal tilings cover the same convex polygon differently (the center parts of the figure).  Instead of gluing one of these tilings to a copy of itself, we may form a different cubic planar graph by adding edges from the corners of one of these zonotopal tilings to the corners of the other, as shown in the bottom right of the figure.  It turns out that this graph is also a cubic partial cube.

A simpler example of the same phenomenon is depicted in Figure~\ref{fig:Q}.  The arrangement in the left of the figure is dual to the zonotopal tiling shown in the center.  If we form a cubic planar graph by connecting corresponding corners of two copies of the tiling, one rotated $90^\circ$ from the other (instead of centrally reflected as our three-dimensional representation of $C_A$ would do), the result is again a partial cube.  We first verified this phenomenon computationally, but a theoretical explanation is hinted at by the right part of the figure: overlaying the zonotopal tiling and its rotation produces another zonotopal tiling, although not a cubic one.

\begin{lemma}
\label{lem:glue-two}
Let $T$ and $T'$ be two zonotopal tilings of the same convex polygon, such that both tilings are dual to simplicial pseudoline arrangements that include the line at infinity,
and such that overlaying the vertices and edges of $T$ and $T'$ produces a third zonotopal tiling $T''$.  Form a planar graph $G$ by connecting corresponding vertices of the outer polygons of $T$ and of the central reflection of $T'$.
Then $G$ is a cubic partial cube.
\end{lemma}

\begin{proof}
$G$ can be formed alternatively as follows: let $A$ be the (non-simplicial) pseudoline arrangement dual to $T''$, and form the (non-cubic) partial cube $C_A$ using the construction in Theorem~\ref{thm:pseudo-pcube}.  We may represent $C_A$ geometrically as above by connecting the corners of $T''$ to those of a reflected copy of itself.  We then replace $T''$ by $T$ on one side of the representation, by removing from each face of $T$ the vertices and edges inside the face that come from $T'$.  Each such replacement preserves the distances among the remaining vertices and edges of the graph, because any path in $T''$ through the interior of a face can be replaced by an equally short path around the boundary of the face.  Similarly, we can replace $T''$ by $T'$ on the other side of the representation, by removing from each face of $T'$ the vertices and edges inside that face coming from $T$.  Thus, $G$ is an isometric cubic subgraph of the partial cube $C_A$, and is therefore itself a cubic partial cube.
\end{proof}

\begin{figure}[t]
\centering
\includegraphics[width=4.5in]{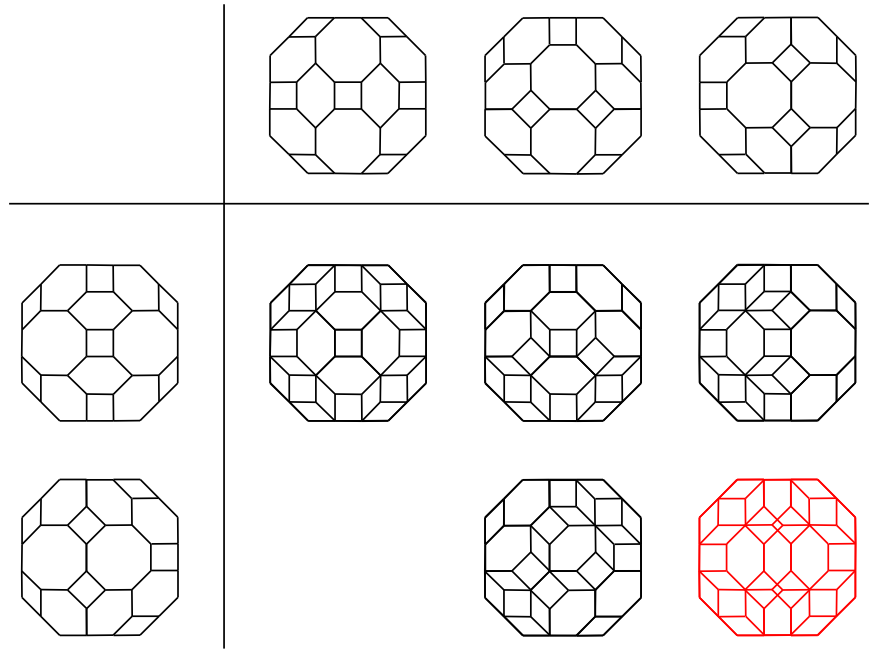}
\caption{Five combinations of two zonotopal tilings of the same polygon and their rotations.  The sixth combination would duplicate one already in the figure, and is omitted.  Four combinations overlay to form zonotopal tilings, and can be connected to form cubic partial cubes by Lemma~\ref{lem:glue-two}.
The fifth combination forms a cubic planar graph that is not a partial cube.}
\label{fig:crosstab}
\end{figure}

Figure~\ref{fig:crosstab} shows all possible combinations of the zonotopal tiling of Figure~\ref{fig:Q} with another zonotopal tiling of the same polygon and with rotations of these tilings.  Four of the combinations yield new cubic partial cubes by Lemma~\ref{lem:glue-two}, but the fifth combination's two tilings do not form a zonotopal tiling when overlaid.  We tested computationally the graph formed by applying the glueing operation of Lemma~\ref{lem:glue-two} to the fifth combination, and found that it is not a partial cube.

Another example of the application of this lemma arises in the arrangement shown in Figure~\ref{fig:zono}.  The zonotopal tiling dual to this arrangement has at its center a decagon, surrounded by a ring of five squares, five hexagons, and five octagons.  This configuration of cells can be rotated by $36^\circ$, forming a second tiling which can be overlayed on the first to produce a finer non-cubic zonotopal tiling.  Connecting the arrangement's original tiling with the modified tiling formed by this rotation produces a non-arrangement-based cubic partial cube with the same number of vertices.

Unfortunately we do not know whether Lemma~\ref{lem:glue-two} may be applied in infinitely many cases, or whether the gluing operation described in the lemma can result in a partial cube even in some cases where the two tilings do not form a zonotopal tiling when overlaid.  The infinite family of pseudoline arrangements described in Figure~\ref{fig:P37} forms zonotopal tilings that may be glued to rotated copies of themselves, forming cubic planar graphs, but these pairs of tilings do not form a single tiling when overlaid so we do not know whether the result of this gluing operation is ever a partial cube.

\section{Conclusions}

We have found a construction of cubic partial cubes from simplicial arrangements, and used this construction to form several new infinite families and sporadic examples of cubic partial cubes.  As a consequence, cubic partial cubes are now much more numerous and varied than previously thought, numerous enough that it is difficult to catalog them exhaustively.

We do not yet know whether there are infinitely many cubic partial cubes that are not of the form $C_A$ for a simplicial line or pseudoline arrangement $A$.  Our Lemma~\ref{lem:glue-two} provides a general framework for the construction of such graphs, but we do not know whether there are infinitely many pairs of zonotopal tilings that can be combined with this lemma.  When it comes to nonplanar graphs, the story is even less clear: it appears that only one nonplanar cubic partial cube is known, having 20 vertices labeled by the set of all five-dimensional bitvectors with two or three bits equal to one.  It would be of interest to determine whether there are any more such examples.

\raggedright
\bibliographystyle{abuser}
\bibliography{simplicial-pcubes}

\end{document}